\documentclass[a4paper,11pt]{article}
\usepackage{amsmath,amssymb,amsthm}
\usepackage{algorithm,algorithmic}

\newcommand{\ie}{\textit{i.e., }}
\renewcommand{\P}{\mathsf{P}}
\newcommand{\E}{\mathsf{E}}
\newcommand{\Var}{\mathsf{Var}}
\newcommand{\R}{\mathbb{R}}
\newcommand{\Ac}{\mathcal{A}}
\newcommand{\Ec}{\mathcal{E}}

\newenvironment{varalgorithm}[1]
  {\algorithm\renewcommand{\thealgorithm}{#1}}
  {\endalgorithm}

\newtheorem{theorem}{Theorem}
\newtheorem{lemma}{Lemma}

\begin{document}
\title{Coloring Random Non-Uniform Bipartite Hypergraphs}
\author{Debarghya Ghoshdastidar, Ambedkar Dukkipati\\
Department of Computer Science \& Automation\\
Indian Institute of Science\\
email: \{debarghya.g,ad\}@csa.iisc.ernet.in
}
\date{\today}

\maketitle

\begin{abstract}
Let $H_{n,(p_m)_{m=2,\ldots,M}}$ be a random non-uniform hypergraph
of dimension $M$ on $2n$ vertices, where the vertices are split into
two disjoint sets of size $n$, and  colored by two distinct colors.  
Each non-monochromatic edge of size $m=2,\ldots,M$ 
is independently added with probability $p_m$.
We show that if $p_2,\ldots,p_M$ are such that the expected number of edges in the
hypergraph is at least $dn\ln n$, for some $d>0$ sufficiently large,
then with probability $(1-o(1))$, one can find a proper 2-coloring of 
$H_{n,(p_m)_{m=2,\ldots,M}}$ in polynomial time.
We present a polynomial time algorithm for hypergraph 2-coloring, and 
provide discussions on extension of the approach for $k$-coloring of non-uniform hypergraphs.
\end{abstract}

\section{Introduction}

A hypergraph $H = (V,E)$ is said to be bipartite or 2-colorable
if the vertex set $V$ can be partitioned into two disjoint sets $V_1$ and $V_2$
such that every edge $e\in E$ has non-empty intersections with both the partitions.
In the case of graphs, one can easily find the two partitions from any given instance of
$H$ by breadth first search.
However, the problem turns out to be notoriously hard if edges of size more than 2 are present.
In fact, in the case of bipartite 3-uniform and 4-uniform hypergraphs,
it is well known that the problem is NP-hard~\cite{Dinur_2005_jour_Combinatorica,Khot_2014_conf_SODA}.

In general, finding a proper 2-coloring is relatively easy if the hypergraph is sparse. 
In an answer to a question asked by Erd\"os~\cite{Erdos_1963_jour_NordikMat} on 2-colorability of uniform hypergraphs, it is now known that for large $m$,
any $m$-uniform hypergraph on $n$ vertices with at most 
$2^m0.7\displaystyle\sqrt{\frac{m}{\ln m}}$ edges is 2-colorable~\cite{Radhakrishnan_1998_conf_FOCS}. As pointed in~\cite{Radhakrishnan_1998_conf_FOCS}, the result can also be extended to 
non-uniform hypergraphs with minimum edge size $m$. 
However, it is much worse if the restriction on the minimum edge size and the
number of hyperedges is not imposed. Even when a hypergraph is 2-colorable, the best
known algorithms~\cite{Alon_1996_jour_NordicJComput,Chen_1996_conf_IPCO}
require $O\left((n\ln n)^{1-1/M}\right)$ colors to properly color the hypergraph
in polynomial time,
where $M$ is the maximum edge size, also called dimension, of the hypegraph.
In recent years, 2-colorability of random hypergraphs has also received considerable attention.
Through a series of works~\cite{Achlioptas_2008_conf_FOCS,CojaOghlan_2012_conf_SODA,Panangiotou_2012_conf_STOC},
it is now established that random uniform hypergraphs are 2-colorable only when
the number of edges are at most $Cn$, for some constant $C>0$.
Thus, it is evident that coloring relatively dense hypergraphs is difficult unless the 
hypergraph admits a ``nice" structure.

In spite of the hardness of the problem,
there are a number of applications that require hypergraph coloring algorithms.
For instance, such algorithms have been used for approximate DNF counting~\cite{Lu_2004_jour_SIAMJDiscMath}, as well as in various resource allocation and scheduling
problems~\cite{Capitanio_1995_jour_IJPP,Ahuja_2002_conf_APPROX}.
 The connection between ``Not-All-Equal" (NAE) SAT
and hypergraph 2-coloring also demonstrate its significance in context of satisfiability problems. 
Among the various approaches studied in the literature, perhaps
the only known non-probabilistic instances of efficient 2-coloring are in the cases 
where the hypergraph is $\alpha$-dense, 3-uniform and bipartite~\cite{Chen_1996_conf_IPCO},  
or where the hypergraph is $m$-uniform and its every edge has equal number of vertices of either colors~\cite{McDiarmid_1993_jour_CombProbComput}.

In this paper, we consider the problem of coloring random non-uniform hypergraphs of dimension $M$,
that has an underlying planted bipartite structure. We present a polynomial time algorithm
that can properly 2-color instances of the random hypergraph with high probability whenever
the expected number of edges in at least $dn\ln n$ for some constant $d>0$.
To the best of our knowledge, such a model has been only considered 
by Chen and Frieze~\cite{Chen_1996_conf_IPCO}, who extended a graph coloring 
approach of Alon and Kahale~\cite{Alon_1997_jour_SIAMJComput} to 
present an algorithm for 
 2-coloring of 3-uniform bipartite hypergraphs with $dn$ number of edges.
To this end, our work generalizes the results  of \cite{Chen_1996_conf_IPCO} to
non-uniform hypergraphs, and it is the first algorithm that is guaranteed to properly color
non-uniform bipartite hypergraphs using only two colors. We also discuss the possible extension 
of our approach to the case of non-uniform $k$-colorable hypergraphs.

\subsection*{The Main Result}
Before stating the main result of this paper, we present the planted model under 
consideration, which 
is based on the model that is studied in~\cite{Ghoshdastidar_2015_arxiv}.
The random hypergraph $H_{n,(p_m)_{m=2,\ldots,M}}$ is generated
on the set of vertices $V = \{1,2,\ldots,2n\}$, which is arbitrarily split into 
two sets, each of size $n$, and the sets are colored with two different colors.
Given a integer $M$, and $p_2,\ldots,p_M\in[0,1]$, the edges of the hypergraph
are randomly added in the following way. All the edges 
of size at most $M$ are added independently, and for any $e\subset V$, 
\begin{align*}
 \P(e\in E) = \left\{ \begin{array}{ll}
     p_m & \text{if } e \text{ is not monochromatic and } {|e|=m}, \\
     0	     & \text{otherwise}.\\
                      \end{array}\right.
\end{align*}
We prove the following result.
%
\begin{theorem}
\label{thm_spec_color}
 Assume $M=O(1)$. There is a constant $d>0$ such that if 
 \begin{equation}
 \sum\limits_{m=2}^M p_m \binom{2n}{m} \geq {dn\ln n}, 
 \end{equation}
 then with probability 
 $(1-o(1))$, Algorithm~\ref{alg} (presented in next section) finds a proper 2-coloring of the random non-uniform bipartite hypergraph $H_{n,(p_m)_{m=2,\ldots,M}}$. 
\end{theorem}
It is easy to see that the expected number of edges in the hypergraph is 
$\Theta\left(\sum_{m=2}^M p_m\binom{2n}{m}\right)$, and so the condition may be stated
in terms of expected number of edges.

\subsection*{Organization of this paper}
The rest of the paper is organized in the following manner.
In Section~\ref{sec_algorithm}, we present our coloring algorithm, followed by a proof of 
Theorem~\ref{thm_spec_color} in Section~\ref{sec_proof}. In the concluding remarks in 
Section~\ref{sec_conclusion}, we provide discussions about the key assumptions made in this work,
and also the possible extensions of our results to $k$-coloring and strong coloring of non-uniform
hypergraphs.  The appendix contains proofs of the lemmas mentioned in Section~\ref{sec_proof}.

\section{Spectral algorithm for hypergraph coloring}
\label{sec_algorithm}

The coloring algorithm, presented below, is similar 
in spirit to the spectral methods of~\cite{Alon_1997_jour_SIAMJComput,Chen_1996_conf_IPCO},
but certain key differences exist, which are essential to deal with
non-uniform hypergraphs. 

Given a hypergraph $H = (V,E)$, 
an initial guess of the color classes is formed by exploiting the spectral properties of a certain matrix
$A\in\mathbb{R}^{|V|\times|V|}$ defined as 
\begin{align}
 A_{ij} = \left\{ \begin{array}{rl}
 \displaystyle\sum_{e\in E: e\ni i,j} \frac{1}{|e|}	& \text{if } i\neq j, \text{ and} \\
 \displaystyle\sum_{e\in E: e\ni i} \frac{1}{|e|}	& \text{if } i= j. 
 \end{array}\right.
 \label{eq_defnA}
\end{align}
The above matrix has been used in the literature to construct the Laplacian of a 
hypergraph~\cite{Bolla_1993_jour_DiscreteMath,Ghoshdastidar_2015_arxiv},
and is also known to be related to the affinity matrix of the star expansion of 
hypergraph~\cite{Agarwal_2006_conf_ICML}. 
The use of matrix $A$ is in contrast to the adjacency based graph construction of~\cite{Chen_1996_conf_IPCO} that is likely to
result in a complete graph if the hypergraph is dense.

The later stage of the algorithm considers an iterative procedure that
is similar 
to~\cite{Alon_1997_jour_SIAMJComput,Chen_1996_conf_IPCO}, but uses a  
weighted summation of neighbors. Such weighting is crucial while
dealing with the edges of
different sizes.

\begin{varalgorithm}{COLOR}
\caption {-- Colors a non-uniform hypergraph $H$:}
\label{alg}
\begin{algorithmic}[1]
 \STATE Define the matrix $A$ as in~\eqref{eq_defnA}.
 \STATE Compute
 $x^A = \underset{\Vert x \Vert_2 = 1}{\textup{arg min~}} x^TAx$.
 \STATE Let $T = \lceil \log_2 n\rceil$, $V_1^{(0)} = \{ i\in V: x_i^A \geq 0\}$ and 
 $V_2^{(0)} = \{ i\in V: x_i^A < 0\}$.
 \FOR {$t = 1,2,\ldots, T$}
 \STATE Let 
 $V_1^{(t)} = \left\{ i\in V: \sum\limits_{j\in V_1^{(t-1)}\backslash\{i\}} A_{ij} <
 \sum\limits_{j\in V_2^{(t-1)}\backslash\{i\}} A_{ij} \right\}$, 
 \newline and $V_2^{(t)} = V\backslash V_1^{(t)}$.
 \ENDFOR
 \IF {{$\exists e\in E$} such that $e\subset V_1^{(T)}$ or $e\subset V_2^{(T)}$}
 \STATE Algorithm FAILS.
 \ELSE
 \STATE 2-Color $V$ according to the partitions 
 $V_1^{(T)},V_2^{(T)}$.
 \ENDIF
\end{algorithmic}
\end{varalgorithm}

\section{Proof of Main Result}
\label{sec_proof}

We now prove Theorem~\ref{thm_spec_color}.
Without loss of generality, assume that the true color classes in $V$ are $\{1,2,\ldots,n\}$ and $\{n+1,\ldots,2n\}$.
Also, let $W^{(t)}$, $t=0,1,\ldots,T$, denote the incorrectly colored
vertices after iteration $t$,
with $W^{(0)}$ being the incorrectly colored nodes after initial spectral step.
We prove Theorem~\ref{thm_spec_color} by showing with probability $(1-o(1))$, 
the size of $W^{(T)} <1$, which implies that all nodes are correctly colored, and hence, the hypergraph must be 
properly colored.

The first lemma bounds the size of $W^{(0)}$, \ie the error incurred at the initial spectral step.
\begin{lemma}
\label{lem_spectral}
With probability $(1-o(1))$,
$|W^{(0)}| \leq \displaystyle\frac{n}{M^22^{2M+4}}$.
\end{lemma}
Next, we analyze the iterative stage of the algorithm to make the following claim,
which characterizes the vertices that are correctly colored after iteration $t$.
\begin{lemma}
\label{lem_iteration_charac}
 Let $\eta = \displaystyle\frac{1}{2^{M+2}}\sum\limits_{m=2}^M\frac{p_m(n-1)}{m}\binom{n-2}{m-2}$. 
 For any $t\in\{1,\ldots,T\}$,
 if $\sum\limits_{j\in W^{(t-1)}\backslash\{i\}} A_{ij} < \eta$ for any $i\in V$, then 
 $P(i\in W^{(t)})\leq n^{-\Omega(d)}$.
\end{lemma}
Note that there are only $T=\lceil \log_2 n\rceil$ iterations, and $|V| =2n$. 
Combining the result of Lemma~\ref{lem_iteration_charac} with union bound, we can conclude
that with probability $(1-o(1))$, for all iterations $t=1,2,\ldots,T$, 
there does not exist any $i\in V$ such that
$\sum\limits_{j\in W^{(t-1)}\backslash\{i\}} A_{ij} < \eta$.
We also make the following observation, where $\eta$ is defined in Lemma~\ref{lem_iteration_charac}.
\begin{lemma}
\label{lem_iteration_size}
With probability $(1-o(1))$, there does not exist $C_1,C_2\subset V$ such that $|C_1|\leq\frac{n}{M^22^{2M+4}}$,
$|C_2| = \frac12 |C_1|$ and for all $i\in C_2$, $\sum\limits_{j\in C_1\backslash\{i\}} A_{ij} \geq \eta$.
\end{lemma}
We now use the above lemmas to proceed with the proof of Theorem~\ref{thm_spec_color}.
Lemma~\ref{lem_spectral} shows that $|W^{(0)}|\leq \frac{n}{M^22^{2M+4}}$ with probability $(1-o(1))$.
Conditioned on this event, and due to the conclusion of Lemma~\ref{lem_iteration_charac},
one can argue that Lemma~\ref{lem_iteration_size} is violated unless
$|W^{(t)}| < \frac12 |W^{(t-1)}|$ for all iteration $t$ with probability $(1-o(1))$. 
Thus, in each iteration,
the number of incorrectly colored vertices are reduced by at least half. Hence, after 
$T=\lceil \log_2 n\rceil$ iterations, $|W^{(T)}| <1$, which implies that all vertices are correctly colored.

\section{Discussions and Concluding remarks}
\label{sec_conclusion}
In this paper, we showed that a random non-uniform bipartite hypergraph of dimension $M$ 
with balanced partitions can be properly 2-colored with 
probability $(1-o(1))$ by a polynomial time algorithm.
The proposed method uses a spectral approach to form initial guess of the color classes,
which is further refined iteratively.
To the best of our knowledge, this is the first work on 2-coloring bipartite non-uniform hypergraphs.
Previous works~\cite{Chen_1996_conf_IPCO,Krivelvich_2003_jour_JAlgo} 
have only restricted to the case of uniform hypergraphs.

\subsection*{A note on the assumptions in Theorem~\ref{thm_spec_color}}

The key assumptions made in this paper are the following:
\begin{enumerate}
\item $M = O(1)$, and 
\item $p_2,\ldots,p_M$ are such that
the expected number of edges is larger than $dn\ln n$, where $d>0$ is a large constant.
\end{enumerate}
The assumption $M = O(1)$ is crucial, particularly in Lemma~\ref{lem_spectral},
and helps to ensure that $d$ can be chosen to be a constant. This can be avoided 
if $d$ is allowed to increase with $n$ appropriately. We note that a previous 
work on spectral hypergraph partitioning~\cite{Ghoshdastidar_2015_arxiv} allows
$M$ to grow with $n$, but imposes an additional restriction so that the number of 
edges of larger size decay rapidly.

The second assumption is stronger than the one in \cite{Chen_1996_conf_IPCO},
where it was shown that a random bipartite 3-uniform hypergraph can be properly
2-colored with high probability if the expected number of edges is $dn$.
This is due to the use of matrix Bernstein inequality~\cite{Tropp_2012_jour_FOCM}
in Lemma~\ref{lem_spectral} that does not provide useful bounds in the most sparse 
case. On the other hand, Chen and Frieze~\cite{Chen_1996_conf_IPCO}
use the techniques of Kahn and Szemeredi~\cite{Friedman_1989_conf_STOC}
that allows them to work in the most sparse regime. 
However, it is not clear how the 
same techniques can be extended even to uniform hypergraphs of higher order.
Thus, it remains an open problem whether a similar result can be proved when the  number of edges in the hypergraph grows linearly with $n$.

\subsection*{$k$-coloring of hypergraphs}
Though Algorithm~\ref{alg} has been presented only for the hypergraph 2-coloring problem,
one may easily extend the approach to achieve a $k$-coloring,
where the objective is to color the vertices of the hypergraph with $k$ colors such that no edge 
is monochromatic.
A possible extension of Algorithm~\ref{alg} is as follows:
\begin{enumerate}
 \item
 In Step~2, compute the eigenvectors corresponding to the $(k-1)$ smallest eigenvalues of $A$. 
 \item
 Use $k$-means algorithm~\cite{Ostrovsky_2013_jour_JACM} to cluster rows of the eigenvector matrix into $k$ groups,
 and define the initial guess for the color classes $V_1^{(0)},\ldots,V_k^{(0)}$ in Step~3 according
 to the above clustering.
 \item
 The iterative computation in Step~6 is modified by defining
 \begin{displaymath}
 \qquad
 V_l^{(t)} = \left\{ i\in V: \sum\limits_{j\in V_l^{(t-1)}\backslash\{i\}} A_{ij} <
 \sum\limits_{j\in V_{l'}^{(t-1)}\backslash\{i\}} A_{ij} \text{ for all } l'\neq l\right\}
 \end{displaymath}
 for $l=1,2,\ldots,(k-1)$, and $V_k^{(t)} = V\backslash \left(\bigcup_{l<k} V_l^{(t)}\right)$.
\end{enumerate}
In the above modification, we borrow the popular idea of using $k$-means on the rows of 
eigenvector matrix to find $k$ planted partitions in a graph or 
hypergraph~\cite{Lei_2015_jour_AnnStat,Ghoshdastidar_2015_arxiv}.

We believe that the result in Theorem~\ref{thm_spec_color} can be extended to this setting,
where the random model allows for $k$ planted color classes in the hypergraph
with non-monochromatic edges generated in the aforementioned manner.
Assuming $k=O(1)$ and $k$-means algorithm always provides a near optimal solution,
one can follow the arguments of~\cite{Ghoshdastidar_2015_arxiv}
to prove a result similar to Lemma~\ref{lem_spectral}. 
On the other hand, Lemmas~\ref{lem_iteration_charac} and~\ref{lem_iteration_size} should hold
for an appropriate choice of $\eta$. Hence, one can comment that the algorithm
achieves a proper $k$-coloring with probability $(1-o(1))$.

We also note that Algorithm~\ref{alg} can be used for finding solutions of NAE-SAT problems.
The extension of~\ref{alg} is also applicable for strong coloring of hypergraphs, which finds 
applications in design of communication networks~\cite{Wu_2015_jour_TIT}.

\appendix
\section*{Proofs of technical lemmas}

\subsection*{Proof of Lemma~\ref{lem_spectral}}
We view the random matrix $A\in\R^{2n\times 2n}$, as a perturbation of its expected value $\Ac=\E[A]$.
Let $\Ec$ denote the collection of all the non-monochromatic subsets of $V$
of size at most $M$.
One can verify that for any $i,j\in V$, $i\neq j$
\begin{align*}
\Ac_{ij} = \sum_{e\in\Ec:e\ni i,j} \frac{p_{|e|}}{|e|}
\qquad\text{and}\qquad
\Ac_{ii} = \sum_{e\in\Ec:e\ni i} \frac{p_{|e|}}{|e|} \;.
\end{align*}
Counting the number of possible edges of each size, one can see that
\begin{align}
\Ac_{ij} = \left\{ \begin{array}{ll}
\alpha_1 - \alpha_2 	& \text{if } i\neq j, \text{ and } i,j \text{ belong to same color class,} \\
\alpha_1			& \text{if } i\neq j, \text{ and } i,j \text{ belong to different color class,} \\
\alpha_1 - \alpha_2 + \alpha_3 & \text{if } i=j,
\end{array}\right.
\end{align}
where
\begin{align*}
\alpha_1 &= \sum_{m=2}^{M}\frac{p_m}{m}\binom{2n-2}{m-2},
\qquad \alpha_2 = \sum_{m=2}^{M}\frac{p_m}{m}\binom{n-2}{m-2},
\\\text{ and } \alpha_3 &= \sum_{m=2}^{M}\frac{p_m}{m}\left(\binom{2n-2}{m-1} - \binom{n-2}{m-1}\right).
\end{align*}
Hence, we can write $\Ac$ as
\begin{align}
 \Ac = \alpha_1 1_{2n\times 2n} - \alpha_2 
 \left(\begin{array}{cc} 1_{n\times n} & 0_{n\times n} \\ 0_{n\times n} & 1_{n\times n} \end{array}\right)
 + \alpha_3 I_{2n},
 \label{eq_Ac}
\end{align}
where $I_{2n}$ is the $2n$-dimensional identity matrix, and $1_{n\times n}$ is a $n\times n$ matrix of all 1's.
One can verify that the smallest eigenvalue of $\Ac$ is $(\alpha_3-n\alpha_2)$, which has multiplicity 1,
and is separated from the other eigenvalues by an eigen-gap of $n\alpha_2$. Moreover, the corresponding 
unit norm eigenvector $x^\Ac$ is such that $x_i^\Ac = \frac{1}{\sqrt{2n}}$ for all $i\leq n$, and 
$x_i^\Ac = -\frac{1}{\sqrt{2n}}$ for all $i> n$, up to a possible change of sign.

At this stage, we refer to a well-known result from matrix perturbation
theory~\cite{Davis_1970_jour_SIAMJNumAnal}. 
We state the result in a particular form that is appropriate in our setting.
The result, as stated in Theorem~\ref{thm_DavisKahan}, has been previously 
used in~\cite[Lemma~4.4]{Ghoshdastidar_2015_arxiv} and~\cite{Lei_2015_jour_AnnStat}.
\begin{theorem}[Davis-Kahan $\sin\Theta$ theorem]
\label{thm_DavisKahan}
Let $\Ac\in\R^{d\times d}$ be a symmetric matrix, and $A$ be an additive perturbation of $\Ac$.
Let $S\subset \R$ be any interval that contains exactly $k$ eigenvalues of $\Ac$.
Define
\begin{displaymath}
\delta = \min\{|\lambda-\lambda'|: \lambda\in S, \lambda'\notin S, \text{ and } 
\lambda,\lambda' \text{ are eigenvalues of } \Ac\}.
\end{displaymath}
If $\delta>2\Vert A-\Ac\Vert_2$, then $S$ also contains exactly $k$ eigenvalues of $A$.

Let $X,\mathcal{X}\in\R^{d\times k}$
be orthonormal eigenvector matrices for the eigenvalues in $S$ of $A,\Ac$ respectively.
Then there is an orthonormal (rotation) matrix $Q\in\R^{k\times k}$ such that
\begin{displaymath}
 \Vert X - \mathcal{X}Q\Vert_F \leq \frac{2\sqrt{2k}\Vert A - \Ac\Vert_2}{\delta} \;.
\end{displaymath}
\end{theorem}

By viewing $A$ as a perturbation of $\Ac$ and noting that the eigen-gap $\delta = n\alpha_2$, 
one can use Theorem~\ref{thm_DavisKahan} 
to conclude that if 
$\alpha_2 > \frac{2}{n} \Vert A - \Ac\Vert_2$, then
\begin{align}
\Vert x^A - x^\Ac\Vert_2 \leq \frac{2\sqrt{2}\Vert A - \Ac\Vert_2}{n\alpha_2}\;.
\label{eq_xA_xAc_bound}
\end{align} 
One can write $A$ as  $A = \sum\limits_{e\in\mathcal{E}} \frac{h_e}{|e|}  a_e a_e^T$, 
where, for each set $e\in\mathcal{E}$, $h_e$ is a Bernoulli$(p_{|e|})$ random variable,
and $a_e\in\{0,1\}^{2n}$ is such that $(a_e)_i = 1$ only when $i\in e$. 
Hence, one may view $A$ as a sum of independent random matrices.
To this end, the following concentration inequality is quite useful to derive a bound
on the perturbation $\Vert A-\Ac\Vert_2$.
\begin{theorem}[Matrix Bernstein inequality~\cite{Tropp_2012_jour_FOCM}]
Consider a finite sequence $X_1,X_2,\ldots,X_L$ of independent, random, self-adjoint matrices with dimension $d$. Assume that each random matrix satisfies $\Vert X_l - \E[X_l]\Vert_2 \leq R$
almost surely. Define $X = \sum\limits_{l=1}^L X_l$, and let $\Var(X) = \E\left[ (X-\E[X])^2\right]$,
where we assume all the above expectations exist. Then for all $t>0$,
\begin{displaymath}
\P\left( \Vert X - \E[X] \Vert_2 \geq t\right) \leq d\exp\left(\frac{-t^2}{2\Var(X) + \frac23Rt}\right).
\end{displaymath}  
\end{theorem}
The above result directly implies
\begin{align}
 \P(\Vert A - \Ac\Vert_2 > 4\sqrt{n\alpha_1 \ln n})
 \leq 4n\exp\left(- \frac{16n\alpha_1 \ln n}{2\Vert \Var(A)\Vert_2 + \frac{8}{3}\sqrt{n\alpha_1 \ln n}}\right).
 \label{eq_AAc_bound}
\end{align}
We note that choosing $d$ large enough, one can satisfy $n\alpha_1 > \ln n$. Also, observe that
\begin{align*}
\Vert \Var(A)\Vert_2  \leq \max_i \sum_{j=1}^{2n} (\Var(A))_{ij} \leq \max_i \sum_{j=1}^{2n} \Ac_{ij} \leq 4n\alpha_1.
\end{align*}
Substituting these in~\eqref{eq_AAc_bound}, we have
\begin{align}
 \P(\Vert A - \Ac\Vert_2 > 4\sqrt{n\alpha_1 \ln n})
 &\leq 4n\exp\left(- \frac{16n\alpha_1 \ln n}{8n\alpha_1 + \frac{8}{3}n\alpha_1}\right)
 \\&= \frac{4}{\sqrt{n}}= o(1).
 \nonumber
\end{align}
Thus, with probability $(1-o(1))$ we have $\Vert A - \Ac\Vert_2 \leq 4\sqrt{n\alpha_1 \ln n}$. 
Due to this bound, one can argue that if $n\alpha_2 > 8\sqrt{\alpha_1 n \ln n}$, \ie
$\frac{\alpha_1}{\alpha_2^2} < \frac{n}{64\ln n}$, then the condition in Theorem~\ref{thm_DavisKahan}
is satisfied, and the preturbation bound~\eqref{eq_xA_xAc_bound} holds.
We can compute that
\begin{align*}
 \frac{\alpha_1}{\alpha_2^2} 
 &= \frac{\sum\limits_{m=2}^M \frac{p_m}{m}\binom{2n-2}{m-2}}{\left(\sum\limits_{m=2}^M \frac{p_m}{m}\binom{n-2}{m-2}\right)^2}
 \\&\leq \frac{n^22^{2M+2}}{\sum\limits_{m=2}^M p_m(m-1)\binom{2n}{m}}
 \\&\leq \frac{n2^{2M+2}}{d\ln n} \;.
\end{align*}
Hence, choosing $d$ sufficiently large, the above mentioned condition holds, and one can claim
from~\eqref{eq_xA_xAc_bound} that
\begin{displaymath}
\Vert x^A - x^\Ac\Vert_2 \leq \displaystyle\frac{8\sqrt{2n\alpha_1 \ln n}}{n\alpha_2}
\leq \frac{2^{M+4.5}}{\sqrt{d}}\;.
\end{displaymath}

Now, we define the set $\widehat{W}\subset V$ as $\widehat{W} = \{i\in V: |x_i^A - x_i^\Ac| \geq \frac{1}{\sqrt{2n}}\}$. 
From the definition of the color classes $V_1^{(0)}, V_2^{(0)}$, it directly follows that any vertex not in $\widehat{W}$
must be correctly colored. Hence, 
\begin{align*}
 |W^{(0)}| &\leq |\widehat{W}| 
 \\&\leq  \sum_{i\in\widehat{W}} 2n |x_i^A - x_i^\Ac|^2
 \\&\leq 2n\Vert x^A - x^\Ac\Vert_2^2 
 \\&= O\left(\frac{n}{d}\right),
\end{align*}
where the bound holds with probability $(1-o(1))$.
Thus, choosing $d$ sufficiently large, one obtains that
$|W^{(0)}| \leq \frac{n}{M^22^{2M+4}}$.

\subsection*{Proof of Lemma~\ref{lem_iteration_charac}}
 Consider any $i\leq n$. Note that $i$ is correctly colored in iteration $t$ if 
 \begin{align*}
  \sum\limits_{j\in V_1^{(t-1)}\backslash\{i\}} A_{ij} < \sum\limits_{j\in V_2^{(t-1)}\backslash\{i\}} A_{ij},
 \end{align*}
 or equivalently,
 \begin{align}
  \sum\limits_{j\in V_1^{(t-1)}\backslash\{i\}} A_{ij} < \frac{1}{2}\sum\limits_{j\neq i} A_{ij}.
  \label{eq_minor_color}
 \end{align}
 Hence, it suffices to show that~\eqref{eq_minor_color} holds under the condition stated in the lemma.
 A similar condition can be stated for $i>n$.
 
We note that $\displaystyle \sum\limits_{j\neq i} A_{ij} =  
 \sum\limits_{e\in\mathcal{E}:e\ni i} h_e\frac{(|e|-1)}{|e|}$, and so, from Bernstein inequality, we have
 \begin{align*}
 \P&\left(\sum_{j\neq i} A_{ij} \leq \left(1-\frac{1}{2^{M+2}}\right) \sum_{j\neq i} \Ac_{ij}\right)
 \\&\leq \exp\left(- \frac{\frac{1}{2^{2M+4}} \left(\sum\limits_{j\neq i} \Ac_{ij}\right)^2}{2\sum\limits_{e\in\mathcal{E}:e\ni i}\frac{(|e|-1)^2}{|e|^2}\Var(h_e) 
 + \frac{2}{3.2^{M+2}}\sum\limits_{j\neq i} \Ac_{ij}}\right)
 \\&\leq \exp\left( - \Omega\left(\sum\limits_{j\neq i} \Ac_{ij}\right)\right)
 \\&\leq n^{-\Omega(d)}.
 \end{align*}
 The second inequality holds since for any $e$, $\frac{(|e|-1)^2}{|e|^2}\Var(h_e) \leq \frac{(|e|-1)}{|e|}\E h_e$,
 and the last inequality is true under the condition of Theorem~\ref{thm_spec_color} since
 \begin{align*}
 \sum_{j\neq i} \Ac_{ij} &= (2n-1)\alpha_1 + (n-1)\alpha_2
 \\&= \sum_{m=2}^M \frac{p_m (m-1)}{2n} \left[ \binom{2n}{m} - 2\binom{n}{m}\right]
 \\&= \Omega(d\ln n).
 \end{align*}
 Denoting  $[n-i] = \{1,\ldots,n\}\backslash i$, \ie the first color class excluding vertex $i$,
 we have $\sum\limits_{j\in [n-i]} A_{ij} = \sum\limits_{e\in\mathcal{E}:e\ni i} h_e \frac{|e\cap [n-i]|}{|e|}$,
 and one can bound
 \begin{align*}
 \P&\left(\sum_{j\in [n-i]} A_{ij} \geq \left(1+\frac{1}{2^{M+2}}\right) \sum\limits_{j\in [n-i]} \Ac_{ij}\right)
 \\&\leq \exp\left(- \frac{\frac{1}{2^{2M+4}} \left(\sum\limits_{j\in [n-i]} \Ac_{ij}\right)^2}{2\sum\limits_{e\in\mathcal{E}:e\ni i}\Var(h_e) \frac{|e\cap U|^2}{|e|^2}
 + \frac{2}{3.2^{M+2}}\sum\limits_{j\in [n-i]} \Ac_{ij}}\right)
\\&\leq n^{-\Omega(d)}.
 \end{align*}
 Thus, with probability $(1-n^{-\Omega(d)})$, we have
 \begin{align*}
 \sum_{j\in [n-i]} A_{ij} &< \left(1+\frac{1}{2^{M+2}}\right) \sum_{j\in [n-i]} \Ac_{ij}
 \\&= \sum_{m=2}^M \frac{p_m(n-1)}{m} \left(1+\frac{1}{2^{M+2}}\right) \left(\binom{2n-2}{m-2}-\binom{n-2}{m-2}\right),
 \end{align*}
 and
 \begin{align*}
 \sum_{j\neq i} A_{ij} &> \left(1-\frac{1}{2^{M+2}}\right) \sum_{j\neq i} \Ac_{ij}
 \\&= \sum_{m=2}^M \frac{p_m}{m}\left(1-\frac{1}{2^{M+2}}\right)\left((2n-1)\binom{2n-2}{m-2}-(n-1)\binom{n-2}{m-2}\right).
 \end{align*}
 Using above relation, we can derive~\eqref{eq_minor_color} since
 \begin{align*}
  \sum_{j\in V_1^{(t-1)}\backslash\{i\}} A_{ij}
  &=  \sum_{j\in W^{(t-1)}\cap V_1^{(t-1)}\backslash\{i\}} A_{ij} + \sum_{j\in V_1^{(t-1)}\backslash(W^{(t-1)}\cap\{i\})} A_{ij} 
  \\&\leq  \sum_{j\in W^{(t-1)}\backslash\{i\}} A_{ij} + \sum_{j\in [n-i]} A_{ij} 
  \\&< \eta + \left(1+\frac{1}{2^{M+2}}\right) \sum_{j\in [n-i]} \Ac_{ij}
  \end{align*}
 The first inequality uses the fact $V_1^{(t-1)}\backslash W^{(t-1)}$ is the set of correctly colored nodes, with true color
 same as $i$. Hence, $V_1^{(t-1)}\backslash(W^{(t-1)}\cap\{i\}) \subset [n-i]$.
 From definition of $\eta$, we have  
  \begin{align*}
  &\sum_{j\in V_1^{(t-1)}\backslash\{i\}} A_{ij}
 \\&\leq \sum_{m=2}^M\frac{p_m(n-1)}{m} \left[ \frac{1}{2^{M+2}}\binom{n-2}{m-2} + 
 \left(1+\frac{1}{2^{M+2}}\right) \left(\binom{2n-2}{m-2}-\binom{n-2}{m-2}\right)\right]
  \\&= \sum_{m=2}^M\frac{p_m(n-1)}{m} \left(1-\frac{1}{2^{M+2}}\right) \left[ \binom{2n-2}{m-2}- \frac12\binom{n-2}{m-2}\right]
  \\&+ \sum_{m=2}^M\frac{p_m(n-1)}{2m}\left[ \frac{1}{2^M}\binom{2n-2}{m-2}- \binom{n-2}{m-2}\right]
  - \sum_{m=2}^M\frac{p_m(n-1)}{m2^{M+3}}\binom{n-2}{m-2}.
 \end{align*}
 One can see that the first term is at most 
 $\frac12 \left(1-\frac{1}{2^{M+2}}\right) \sum_{j\neq i} \Ac_{ij}
 < \frac12  \sum_{j\neq i} A_{ij}$.
 On the other hand, we note that 
 \begin{align*}
  \frac{\binom{2n-2}{m-2}}{\binom{n-2}{m-2}} 
  \leq \frac{1}{4} \frac{\binom{2n}{m}}{\binom{n}{m}} 
  \leq \frac14 \frac{\frac{(2n)^m}{m!}}{\frac{n^m}{4.m!}} = 2^m \leq 2^M\;.
 \end{align*}
 So the second term is negative, which proves~\eqref{eq_minor_color}, and the claim follows.

\subsection*{Proof of Lemma~\ref{lem_iteration_size}}
Let $C_1,C_2\subset V$ be arbitrary such that $|C_2| = b$, 
and $E_{C_1C_2}$ be the set of all non-monochromatic subsets of $V$ 
of size at most $M$
that have non-empty intersection with both $C_1$ and $C_2$. 
Then
\begin{align*}
\sum_{e\in E_{C_1C_2}} h_e 
&\geq \frac1M \sum_{e\in E_{C_1C_2}} h_e \frac{|e\cap C_1||e\cap C_2|}{|e|}
\\&\geq \frac1M \sum_{i\in C_2}\sum_{j\in C_1\backslash\{i\}} A_{ij} \geq \frac{b\eta}{M},
\end{align*}
where the last inequality holds under the condition stated in the lemma. Now we bound the probability
\begin{align}
&\P\left(\exists C_1,C_2 \subset V, |C_2|=\frac12 |C_1| \leq \frac{n}{M^22^{2M+5}}, 
\sum_{j\in C_1\backslash\{i\}} A_{ij} \geq \eta ~\forall i\in C_2\right)
\label{eq_A_ij_prob_bound1}
\\&\leq \sum_{b=1}^{\frac{n}{M^22^{2M+5}}}
\P\left(\exists C_1,C_2 \subset V, |C_2|=\frac12 |C_1| = b, \text{ and }
\sum_{e\in E_{C_1C_2}} h_e \geq \frac{b\eta}{M} \right)
\nonumber
\\&\leq \sum_{b=1}^{\frac{n}{M^22^{2M+5}}} \sum_{C_2: |C_2|=b} \sum_{C_1:|C_1|=2b} 
\P\left(\sum_{e\in E_{C_1C_2}} h_e \geq \frac{b\eta}{M} \right)
\nonumber
\end{align}
We observe that
\begin{align*}
 \sum_{e\in E_{C_1C_2}} \E [h_e]
 &=\sum_{m=2}^M \sum_{e\in E_{C_1C_2},|e|=m} p_m
 \\&\leq 2b^2 \sum_{m=2}^M p_m \binom{2n-2}{m-2}
 \\&\leq b^2 2^{M+1}\sum_{m=2}^M p_m \binom{n-2}{m-2}
 \\&\leq \frac{b^2\eta M 2^{2M+4}}{n},
\end{align*}
and the above bound is smaller than $\frac{b\eta}{2M}$ for 
$b\leq\frac{n}{M^22^{2M+5}}$.
Hence, we can write
\begin{align*}
&\P\left(\sum_{e\in E_{C_1C_2}} h_e \geq \frac{b\eta}{M} \right)
\\& \leq\exp\left(\frac{-\left(\frac{b\eta}{M} -\sum_{e\in E_{C_1C_2}}\E[h_e]\right)^2}{
2\sum_{e\in E_{C_1C_2}} \Var(h_e) + \frac23 \left(\frac{b\eta}{M} - 
\sum_{e\in E_{C_1C_2}}\E[h_e]\right)}\right)
\\&\leq \exp\left(-\frac{3b\eta}{16M}\right).
\end{align*}
Substituting in~\eqref{eq_A_ij_prob_bound1}, we have the probability of the existence of
$C_1,C_2$ with mentioned conditions is at most
\begin{align*}
\sum_{b=1}^{\frac{n}{M^22^{2M+5}}} \binom{2n}{b}\binom{2n}{2b}\exp\left(-\frac{3b\eta}{16M}\right)
&\leq \sum_{b=1}^{\infty} \left(2n\exp\left(1-\frac{\eta}{16M}\right)\right)^{3b}.
\end{align*}
Under the assumption of Theorem~\ref{thm_spec_color}, one can verify that
$\eta\geq \frac{d\ln n}{2^{2M+4}}$.
So for large $d$, the above geometric series converges, and is at most $n^{-\Omega(d)} = o(1)$. Hence, the claim.

\end{document}